\documentstyle[12pt,xypic]{article}
\def\Im{\mathop{\rm Im}\nolimits}
\def\Re{\mathop{\rm Re}\nolimits}
\def\Sp{\mathop{\rm Sp}\nolimits}
\def\tr{\mathop{\rm tr}\nolimits}

\def\e{{\varepsilon}}
\def\ep{{\epsilon}}

\def\b{{\beta}}
\def\s{{\sigma}}
\def\la{{\langle}}
\def\ra{{\rangle}}

\def\arr{{\longrightarrow}}

\def\d{{\delta}}
\def\g{{\gamma}}
\def\l{{\lambda}}
\def\ph{{\varphi}}
\def\o{{\omega}}
\def\i{{\infty}}

\def\G{{\Gamma}}
\def\q{{\quad$\Box$\par\medskip}}
\def\be{\begin{equation}}
\def\ee{\end{equation}}

\newtheorem{thm}{Theorem}
\newtheorem{lem}[thm]{Lemma}
\newtheorem{prop}[thm]{Proposition}
\newtheorem{dfn}[thm]{Definition}
\newtheorem{cor}[thm]{Corollary}

\newcommand{\norm}[1]{\left\| #1 \right\|}
\newcommand{\nnn}[1]{|\! |\! | #1 |\! |\! |}
\newcommand{\ov}[1]{\overline{#1}}

\date{}
\author{V.~M.~Manuilov}
\title{Almost representations and asymptotic representations of discrete
groups}
\begin{document}
\maketitle

 \begin{abstract}
We define for discrete finitely presented groups a new property related to
their asymptotic representations. Namely we say that a groups has the
property AGA if every almost representation generates an asymptotic
representation. We give examples of groups with and without this property.
For our example of a group $\G$ without AGA the group
$K^0(B\G)$ cannot be covered by asymptotic representations of $\G$.

 \end{abstract}

One of the reasons of attention to almost and asymptotic representations
of discrete groups \cite{Harp-Kar,c-hig} is their relation to $K$-theory
of classifying spaces \cite{c-hig,mish-noor}. It was shown in
\cite{mish-noor} that in the case of finite-dimensional classifying space
$B\G$ in order to construct a vector bundle over it out of an asymptotic
representation of $\G$ it is sufficient to have an $\e$-almost
representation of $\G$ with small enough $\e$. Of course an $\e$-almost
representation contains less information than the whole asymptotic
representation, but it turns out that often the information contained in
an $\e$-almost representation makes it possible to construct the
corresponding asymptotic representation. In the present paper we give the
definition of this property, prove this property for some classes of
groups and finally give an example of a group without this property. We
discuss also this example in relation with its $K$-theory.

\section{Basic definitions}

Let $\G$ be a finitely presented discrete group, and let
$\G=\la F|R\ra=\la g_1,\ldots,g_n|r_1,\ldots,r_k\ra$ be its presentation
with $g_i$ being generators and $r_j=r_j(g_1,\ldots,g_n)$ being relations.
We assume that the set $F=\{g_1,\ldots,g_n\}$ is symmetric, i.e. for every
$g_i$ it contains $g_i^{-1}$ too, and the set $R$ of relations contains
relations of the form $g_ig_i^{-1}$, though we usually will skip these
additional generators and relations.

By $U_\i$ we denote the direct limit of the groups $U_n$
with respect to the natural inclusion $U_n\arr U_{n+1}$ supplied with the
standard operator norm. The unit matrix we denote by $I\in U_\i$.

 \begin{dfn}
{\rm
A set of unitaries $\s=\{u_1,\ldots,u_n\}\subset U_\i$ is called an {\it
$\e$-almost representation} of the group $\G$ if after substitution of
$u_i$ istead of $g_i$, $i=1,\ldots,n$, into $r_j$ one has
 $$
\norm{r_j(u_1,\ldots,u_n)-I}\leq\e
 $$
for all $j=1,\ldots,k$.

}
 \end{dfn}

In this case we write $\s(g_i)=u_i$. Remark that this definition depends
on a choice of presentation of the group $\G$, but we will see that this
dependence is not important. Let $\la h_1,\ldots,h_m|s_1,\ldots,s_l\ra$ be
another presentation of $\G$. For an $\e$-almost representation $\s$ with
respect to the first presentation we can define the set of unitaries
$v_1\ldots,v_m\in U_\i$, $v_i=\ov{\s}(h_i)$ putting
$\ov{\s}(h_i)=\s(g_{j_1})\cdot\ldots\cdot\s(g_{j_{n_i}})$,
where $h_i=g_{j_1}\cdot\ldots\cdot g_{j_{n_i}}$.
By the same way starting
from the set $\ov{\s}(h_i)$ we can construct the set $\ov{\ov{\s}}(g_i)$.

 \begin{lem}\label{nezavis}
There exist constants $C$ and $D$ (depending on the two presentations)
such that $\ov{\s}$ is a $C\e$-almost representation with respect to the
second presentation of $\G$ and for all $g_i$, $i=1,\ldots,n$, one has
$\norm{\ov{\ov{\s}}(g_i)-\s(g_i)}\leq D\e$.

 \end{lem}

{\bf Proof.}
We have to estimate the norms $\norm{s_q(v_1,\ldots,v_m)-I}$,
$q=1,\ldots,l$. To do so notice that every relation $s_q$ can be written
in the form
 \begin{equation}\label{kommut}
s_q=a_1^{-1}r_{j_1}^{\ep_1}a_1\cdot\ldots\cdot
a_{m_q}^{-1}r_{j_{m_q}}^{\ep_{m_q}}a_{m_q}
 \end{equation}
for some $a_i\in \G$, where $\ep_i=\pm 1$.
Let $M$ be the maximal length of the words
$a_i=a_i(g_1,\ldots,g_n)$. Put $b_i'=a_i^{-1}(g_1,\ldots,g_n)\in U_\i$,
$b_i=a_i(g_1,\ldots,g_n)\in U_\i$. Then one has
 $$
\norm{b_i'b_i-I}\leq M\e.
 $$
It follows from (\ref{kommut}) that
 $$
s_q(v_1,\ldots,v_m)=
b_1'r_{j_1}(u_1,\ldots,u_n)b_1\cdot\ldots\cdot
b_{m_q}'r_{j_{m_q}}(u_1,\ldots,u_n)b_{m_q},
 $$
but as for every $i$ one has
 $$
\norm{b_i'r_{j_i}^{\ep_i}(u_1,\ldots,u_n)b_i-I}\leq
\norm{b_i'b_i-I}+\norm{r_{j_i}^{\ep_i}(u_1,\ldots,u_n)-I}\leq(M+1)\e,
 $$
so
 $$
\norm{s_q(v_1,\ldots,v_m)-I}\leq m_q(M+1)\e,
 $$
which proves the first statement of the Lemma. The second statement is
proved in a similar way. \q

As the number of generators is finite, so the image of every almost
representation lies in finite matrices, $\s\in U_n$ for some $n$. The
minimal such $n$ is called a dimension of $\s$. Usually we will ignore the
remaining infinite unital tail of the matrices $\s(g_i)$ and write
$\s(g_i)\in U_n$ instead of $U_\i$.

The set of all $\e$-almost representations of the group $\G$ we denote by
$R_\e(\G)$.
The deviation of an almost representation can be measured by the
value
 $$
\nnn{\s}=\max_j\norm{r_j(u_1,\ldots,u_n)-I}.
 $$
Notice that both these definitions also depend on the choice of a
presentation of $\G$.

\medskip
 \begin{dfn}[{\rm cf. \cite{c-hig}}]
{\rm
A set of norm-continuous unitary paths
$\s_t=\{u_1(t),\ldots,u_n(t)\}\subset U_\i$, $t\in[0,\i)$, is called an
{\it asymptotic representation} of the group $\G$ if $\nnn{\s_t}$ tends to
zero when $t\to\i$.

}
 \end{dfn}

Due to the Lemma \ref{nezavis} this definition does not depend on the
choice of a presentation of $\G$. The set of all asymptotic representations
of the group $\G$ we denote by $R_{asym}(\G)$.

\medskip
Now we are ready to define a new property of finitely generated groups
which we call AGA ({\it Almost} representations {\it Generate Asymptotic}
representations).

 \begin{dfn}
{\rm
A group $\G$ possesses the property AGA if for every $\e>0$
one can find a number $\d(\e)$ (with the property $\d(\e)\to 0$
when $\e\to 0$) such that for every
almost representation $\s\in R_\e(\G)$  there exists an asymptotic
representation $\s_t\in R_{asym}(\G)$ such that $\s_0=\s$ and
$\nnn{\s_t}\leq \d(\e)$ for all $t\in [0,\i)$.

}
\end{dfn}

 \begin{lem}
The property AGA does not depend on the choice of a presentation of the
group.

 \end{lem}

{\bf Proof } immediately follows from the Lemma \ref{nezavis}. \q

 \begin{thm}
The following groups have the property AGA:
 \begin{enumerate}
\vspace{-\itemsep}
\item
free groups,

\vspace{-\itemsep}
\item
free products of groups having the property AGA,

\vspace{-\itemsep}
\item
subgroups of finite index in groups having the property AGA,

\vspace{-\itemsep}
\item
finite groups,

\vspace{-\itemsep}
\item
free abelian groups,

\vspace{-\itemsep}
\item
fundamental groups of two-dimensional oriented manifolds.

\vspace{-\itemsep}

 \end{enumerate}
 \end{thm}

{\bf Proof.} The first item is obvious --- free groups have no relations,
so every almost representation is a genuine representation.
The same argument works for the second item too. The third item can be
proved in the same way as the Lemma \ref{nezavis} as subgroups of finite
index are defined by a finite number of relations originated from
relations of the group.
The fourth item was proved in \cite{Harp-Kar} --- for finite groups there
exists a genuine representation close to every almost representation.
The fifth item is non-trivial too and was proved in \cite{manCopen}. It
follows also from the Super Homotopy Lemma of \cite{BEEK}.
The sixth item we prove in the next section.

\section{Case of fundamental groups of oriented surfaces}

Let $\G=\la a_1,b_1,\ldots,a_m,b_m|a_1b_1a_1^{-1}b_1^{-1}\cdot\ldots\cdot
a_mb_ma_m^{-1}b_m^{-1}\ra$. Let $\s$ be an $\e$-almost representation of
$\G$, $u_i=\s(a_i)$, $v_i=\s(b_i)$, $u_i,v_i\in U_n$, $i=1,\ldots,m$.
Denote $\g(u,v)=uvu^{-1}v^{-1}$, then we have
 $$
\norm{\g(u_1,v_1)\cdot\ldots\cdot \g(u_m,v_m)-I}\leq\e.
 $$

Consider the map
 \begin{equation}\label{map_gamma}
\g:U_n\times U_n\arr SU_n.
 \end{equation}

To prove the property AGA for fundamental groups of oriented
two-dimensional manifolds we have to use the following elementary
statement about the map (\ref{map_gamma}).

 \begin{lem}\label{lem:podniat}
Let $(u_0,v_0)\in U_n\times U_n$ and let $c(t)\in SU_n$, $t\in [0,1]$, be
a path such that $\g(u_0,v_0)=c(0)$. Then for any $\d>0$ there
exists a path $(u_t,v_t)\in U_n\times U_n$ such that
$\norm{\g(u_t,v_t)-c(t)}<\d$.

 \end{lem}

{\bf Proof.}
Remember that a pair $(u,v)\in U_n\times U_n$ is called {\it irreducible}
if there is no common invariant subspace for $u$ and $v$.
It was shown in \cite{exel} that the set of regular points for the map
$\g$ (\ref{map_gamma}) coincides with the set of irreducible pairs.
Denote the set of reducible pairs $(u,v)\in U_n\times U_n$ by $S$.
For any $k=1,\ldots,n-1$ by $\Sigma_k\subset SU_n$ denote the set of
block-diagonal matrices
$\left(\begin{array}{cc}c_1&0\\0&c_2\end{array}\right)$ with respect to
some invariant subspace $V$, $\dim V=k$, such that $c_1\in SU_k$, $c_2\in
SU_{n-k}$. Put $\Sigma=\cup_k\Sigma_k\subset SU_n$. Then obviously
$\g(S)\subset \Sigma$. Notice that every $\Sigma_k$ is a submanifold
in $SU_n$ with codimension one. So $\Sigma$ divides $SU_n$ into a finite
set of closed path components $M_j$, $\cup_j M_j=SU_n$ and for every point
$c\in M_j$ the set $\g^{-1}(c)$ consists only of regular points. Hence
every path in $M_j$ transversal to its boundary can be lifted up to a path
in $U_n\times U_n$ with a fixed starting point.

Without loss of generality we can assume that the path $c(t)$ is
transversal to every $\Sigma_k$. Let $t_0\in \{c(t)\}\cap\Sigma_k$. It
remains to show that we can lift the path $c(t)$ in some neighborhood of
the point $c_0=c(t_0)$. Let $(u_0,v_0)\in U_n\times U_n$ be such point
that $\g(u_0,v_0)=c_0$. If the point $(u_0,v_0)$ is a regular point then
the statement is obvious. Otherwise we can write
 $$
u_0=\left(\begin{array}{cc}u_1&0\\0&u_2\end{array}\right),\qquad
v_0=\left(\begin{array}{cc}v_1&0\\0&v_2\end{array}\right)
 $$
with respect to some basis and we can assume that matrices $v_1$ and $v_2$
are diagonal. Let
 $$
e^{2\pi i\ph_1},\ldots,e^{2\pi i\ph_k}\quad{\rm and}\quad
e^{2\pi i\ph_{k+1}},\ldots,e^{2\pi i\ph_n}
 $$
be the eigenvalues of
$v_1$ and $v_2$ respectively.
Slightly changing $v_0$ we can assume that for all $i=2,\ldots,k$,
$j=k+2,\ldots,n$ the values $\ph_i-\ph_1$, $\ph_j-\ph_{k+1}$ differ from
each other. Multiplying $v_1$ by $e^{-2\pi i\ph_1 t}$ and $v_2$ by
$e^{-2\pi i\ph_{k+1} t}$, $t\in [0,1]$, we connect the matrix $v_0$ with
the matrix
 $$
v'_0=\left(\begin{array}{cc}v'_1&0\\0&v'_2\end{array}\right)=
\left(\begin{array}{cc}e^{-2\pi i\ph_1}v_1&0\\
0&e^{-2\pi i\ph_{k+1}}v_2\end{array}\right)
 $$
which has two eigenvalues equal to one and all other eigenvalues being
different from each other. Obviously the value
$\g(u_0,v_0)=\g(u_0,v'_0)=c_0$ does not change along this path.
Denote by $e_1,\ldots,e_n$ the basis consisting of the eigenvalues of
$v'_0$ and let $w(t)\in U_n$, $t\in[0,1]$, be a rotation of the vectors
$e_1$ and $e_{k+1}$:
 $$
w(t)e_1=\cos t e_1-\sin t e_{k+1},\quad
w(t)e_{k+1}=\sin t e_1+\cos t e_{k+1},\quad
w(t)e_j=e_j \quad {\rm for}\quad j\neq 1,k+1.
 $$
Obviously $w(t)$ commutes with $v'_0$. Put $u_t=u_0w(t)$. Then
 $$
\g(u_t,v'_0)=u_0w(t)v'_0w^{-1}(t)u_0^{-1}(v'_0)^{-1}=\g(u_0,v'_0)=c_0
 $$
and for $\sin t\neq 0$ the pair $(u_t,v'_0)$ is {\it irreducible}
(since $v'_0$ is diagonal with
only two coinciding eigenvalues, so its invariant subspaces are easy to
describe, then it is easy to check that they are not invariant under the
action of $u_t$) with the same value of $\g$. Then it is possible to
extend the path $(u_t,v_t)$ through the point $c_0$. \q

 \begin{prop}\label{homot_to_1}
Any $\e$-almost representation of $\G$ is homotopically equivalent in
$R_{2\e}(\G)$ to an $2\e$-almost representation with
$u_2=v_2=\ldots=u_m=v_m=I$.

 \end{prop}

{\bf Proof.}
Connect the matrix $\g(u_m,v_m)$ with $I$ by a path $c_m(t)$. Then by the
Lemma \ref{lem:podniat} we can find a path $(u_m(t),v_m(t))\in U_n\times
U_n$ such that
 $$
\norm{\g(u_m(t),v_m(t))-c_m(t)}\leq\frac{\e}{2m}.
 $$
Notice
that the set $\g^{-1}(I)=\{(u,v):uv=vu\}$ is path-connected, so we can
assume that the end point of the path $(u_m(t),v_m(t))$ is $(I,I)$.
Put
 $$
c_{m-1}(t)=\g(u_{m-1},v_{m-1})c^{-1}_m(t).
 $$
Again by the Lemma
\ref{lem:podniat} we can find a path $(u_{m-1}(t),v_{m-1}(t))\in
U_n\times U_n$ such that
 $$
\norm{\g(u_{m-1}(t),v_{m-1}(t))-c_{m-1}(t)}\leq\frac{\e}{2m}.
 $$
Then
 $$
\norm{\g(u_1,v_1)\cdot\ldots\cdot\g(u_{m-1}(t),v_{m-1}(t))
\cdot\g(u_m(t),v_m(t))-I}\leq\e+\frac{\e}{m}
 $$
and at the end point we have $(u_m(t),v_m(t))=(I,I)$. Proceeding by
induction we finish the proof.\q

It now follows from the proposition \ref{homot_to_1} that the property AGA
for the group $\G$ follows from the same property for the group ${\bf Z}^2$
with generators $u_1,v_1$. \q

\section{Example of a group without AGA}

Consider the group $\G=\la a,b,c|aca^{-1}c^{-1},b^2,(ab)^2\ra$.

 \begin{thm}
The group $\G$ does not possess the property AGA.

 \end{thm}

{\bf Proof.}
Put $\o=e^{2\pi i/n}$ and define a
family of almost representations $\s_n$ taking values in $U_n$ by
 $$
\s_n(a){=}\!\left(\begin{array}{cccccc}
\o&&&&&\\
&\o^2&&&&\\
&&\cdot&&&\\
&&&\cdot&&\\
&&&&\cdot&\\
&&&&&\o^n
\end{array}\right)\!\!,\
\s_n(c){=}\!\left(\begin{array}{cccccc}
0&&&&&1\\
1&0&&&&\\
&1&\cdot&&&\\
&&\cdot&\cdot&&\\
&&&\cdot&\cdot&\\
&&&&1&0
\end{array}\right)\!\!,\
\s_n(b){=}\!\left(\begin{array}{cccccc}
&&&&&1\\
&&&&1&\\
&&&\cdot&&\\
&&\cdot&&&\\
&\cdot&&&&\\
1&&&&&
\end{array}\right)\!\!.
 $$

Here the matrices $\s_n(a)$ and $\s_n(c)$ are the Voiculescu matrices
\cite{voi} with the winding number \cite{e-l} equal to one, and one has
 $$
\s_n(a)\s_n(c)\s_n(a)^{-1}\s_n(c)^{-1}=\o\cdot I,\quad
\s_n(b)^2=I,\quad (\s_n(a)\s_n(b))^2=\o\cdot I,
 $$
so
 $$
\e_n=\nnn{\s_n}=|\o-1|\to 0 \quad{\rm when}\ \ n\to\i,
 $$
hence for every $\e>0$ there exists an $\e$-almost representation $\s$ of
$\G$ such that the winding number of the pair $(\s(a),\s(c))$ equals one.

\smallskip
Suppose the opposite, i.e. that the group $\G$ has AGA. Then there should
be such $\e$ that the number $\d(\e)\leq 1$. Take such $\e$ and an
$\e$-almost representation $\s_0$ with a non-zero winding number of the
pair $(\s_0(a),\s_0(c))$.

\smallskip
By supposition there exists an asymptotic representation
$\s_t\in R_{asym}(\G)$ extending $\s_0$ such that
$\nnn{\s_t}\leq 1$ for all $t\in [0,\i)$ and for any $\e'>0$ there exists
$t_0$ such that $\nnn{\s_{t_0}}\leq \e'$. Fix this $t_0$ and denote
$\s_{t_0}$ by $\s$. Let $n$ and $n+m$ be the dimension of the almost
representation $\s_0$ and $\s$ respectively.
Then one has
 $$
\norm{\s(a)\s(c)-\s(c)\s(a)}\leq\e',
\quad
\norm{\s(b)^2-I}\leq\e',
\quad
\norm{\s(a)\s(b)\s(a)\s(b)-I}\leq\e'.
 $$
Notice that as along the whole path $\s_t$ one has
 $$
\norm{\s_t(b)^2-I}\leq\nnn{\s_t}\leq\d(\e)\leq 1,
 $$
so the eigenvalues of $\s_t(b)$ satisfy the estimate $|\l^2-1|\leq 1$,
hence the number of eigenvalues $\l\in\Sp \s_t(b)$ with $\Re\l<0$
does not change along the whole path $\s_t$ in $U_\i$,
therefore the number of eigenvalues $\l\in\Sp\s(b)$ with $|\l+1|\leq \e'$
cannot exceed $n$ (the maximal number of eigenvalues with $\Re\l<0$ of
$\s_0(b)$) and the number of eigenvalues with $|\l-1|\leq \e'$
is not less than $m$.
Then there exists a matrix $\s(b)'\in U_{n+m}$ such that
$\norm{\s(b)-\s(b)'}\leq\e'$ and that the matrix $\s(b)'$ has not more
than $n$ eigenvalues equal to $-1$ and not less than $m$ eigenvalues
equal to $1$. Then we have
 \begin{equation}\label{1*}
\norm{\s(a)\s(b)'\s(a)\s(b)'-I}\leq 3\e'.
 \end{equation}
Notice that $(\s(b)')^2=I$ and
 \begin{equation}\label{1a*}
|\tr(\s(b)')|\geq m-n.
 \end{equation}

Let
 $$
\s(a)=\left(\begin{array}{ccc}
\o_1&&\\
&\ddots&\\
&&\o_{n+m}
\end{array}\right)
 $$
be the matrix of the operator $\s(a)$ in the basis consisting of its
eigenvectors. It was shown in \cite{manFA,BEEK} that if the winding
number of the pair $(\s(a),\s(c))$ is non-zero then for any $\e'>0$ there
exists $\d'(\s')$ such that $\d'(\e')\to 0$ when $\e'\to 0$ and that all
lacunae in $\Sp\s(a)$ do not exceed $\d'(\e')$.
Denote the number of eigenvalues of $\s(a)$ with $|\Im\o_j|>2\e'$ by $N$.
Then we have
 \begin{equation}\label{1d*}
N\geq\frac{2\pi-10\e'}{\d'(\e')}.
 \end{equation}

As $(\s(b)')^2=I$, so it follows from (\ref{1*}) that
 \begin{equation}\label{1c*}
\norm{\s(a)\s(b)'-\s(b)'\s(a)^*}\leq 3\e'.
 \end{equation}
Denote by $b_{ij}$ the matrix elements of the matrix $\s(b)'$.
It follows from (\ref{1c*}) that all matrix elements of
$\s(a)\s(b)'-\s(b)'\s(a)^*$ do not exceed $3\e'$, i.e.
 \begin{equation}\label{2*}
|b_{ii}(\o_i-\ov{\o}_i)|\leq 3\e', \quad i=1,\ldots,n+m.
 \end{equation}
Let us estimate $\tr(\s(b)')$.
We have
 $$
|\tr(\s(b)')|=\left|\sum_{i=1}^{n+m}b_{ii}\right|\leq
\sum_{i=1}^{n+m}|b_{ii}|=\sum\nolimits'|b_{ii}+\sum\nolimits''|b_{ii}|,
 $$
where $\sum'$ denotes the sum for those numbers $i$ for which one has
$|\Im\o_i|>2\e'$ and $\sum''$ is the sum for the remaining numbers.
As for all $i$ one has $|b_{ii}|\leq 1$, so the last sum do not exceed
the number of summands,
 $$
\sum\nolimits''|b_{ii}|\leq n+m-N.
 $$
It follows from (\ref{2*}) that for those $i$ which are included into the
first sum we have $|\o_i-\ov{\o}_i|>4\e'$, hence those $b_{ii}$ satisfy
 $$
|b_{ii}|<\frac{3}{4},
 $$
so
 $$
\sum\nolimits'|b_{ii}|<\frac{3}{4}N,
 $$
hence we have
 $$
|\tr(\s(b)')|<\frac{3}{4}N+n+m-N=n+m-\frac{N}{4}
 $$
and it follows from (\ref{1d*}) that
 \begin{equation}\label{2a*}
|\tr(\s(b)')|<n+m-\frac{\pi-5\e'}{2\d'(\e')}.
 \end{equation}
If we take $\e'$ small enough then $\d'(\e')$ is small enough too and we
get ${\displaystyle\frac{\pi-5\e'}{2\d'(\e')}>2n}$, then (\ref{1d*}) and
(\ref{2a*}) give a contradiction. \q

 \begin{cor}\label{RG=RH}
Let $\s_t\in R_{asym}(\G)$ be an asymptotic representation. Then the
winding number of the pair $(\s_t(a),\s_t(c))$ is zero for big enough $t$.
In particular, it means that $\s_t$ is homotopic to an asymptotic
representation $\rho_t\in R_{asym}(\G)$ with $\rho_t(c)=I$ in the class of
asymptotic representations. \q

\end{cor}

Denote the Grothendieck group of homotopy classes of asymptotic
representations of the group $\G$ by ${\cal R}_{asym}(\G)$. Let $H$ denote
the subgroup $\la a,b|b^2,(ab)^2\ra\cong {\bf Z}_2\ast{\bf Z}_2$. Then
$\G\cong {\bf Z}^2\ast_{\bf Z}H$, so we have
$B\G={\bf T}^2\cup BH$, $S^1={\bf T}^2\cap BH$,
where $B\G$, $BH$, $S^1$ and ${\bf T}^2$ are the classifying spaces of the
groups $\G$, $H$, ${\bf Z}$ and ${\bf Z}^2$ respectively, and the inclusion
$S^1\subset {\bf T}^2$ is the standard inclusion onto the first coordinate.
Then one has an exact sequence
 \begin{equation}\label{K6}
\diagram
K^0(B\G)\rto & K^0(BH)\oplus K^0({\bf T}^2)\rto & K^0(S^1)\dto\\
K^1(S^1)\uto & K^1(BH)\oplus K^1({\bf T}^2)\lto & K^1(B\G)\lto
\enddiagram
 \end{equation}
and as the maps $K^*({\bf T}^2)\arr K^*(S^1)$ are onto, so the vertical
maps in (\ref{K6}) are zero and the group $K^0(B\G)$ contains an element
$\beta$ which is mapped onto the Bott generator of $K^0({\bf T}^2)$.

Remember that in \cite{mish-noor} a map
 \begin{equation}\label{mapasym}
{\cal R}_{asym}(\G)\arr K^0(B\G)
 \end{equation}
was constructed, which factorizes \cite{man-mish} through the assembly map
 $$
{\cal R}_{\cal Q}(\G\times{\bf Z})\arr K^0(B\G),
 $$
where ${\cal R}_{\cal Q}(\G)$ denotes the Grothendieck group of
representations of $\G$ into the Calkin algebra ${\cal Q}$.

It follows from the Corollary \ref{RG=RH} that ${\cal R}_{asym}(\G)={\cal
R}_{asym}(H)$, therefore the element $\beta\in K^0(B\G)$ does not lie in
the image of the map (\ref{mapasym}), hence we obtain

 \begin{cor}
The map (\ref{mapasym}) is not a rational epimorphism for the group $\G$.
\q

 \end{cor}

On the other hand we should remark that the element $\b\in K^0(B\G)$ can be
obtained as an image of a representation of $\G\times{\bf Z}$ into the
Calkin algebra. To get such representation we can take $\s(b)$ to be the
infinite direct sum of matrices of the form
$\left(\begin{array}{cc}0&1\\1&0\end{array}\right)$ and to rearrange the
basis for the matrices $\s_n(a)$ and $\s_n(c)$ in such a way that $\s_n(a)$
and $\s(b)$ would almost commute. We should also insert a number of
intermediate matrices between $\s_n$ and $\s_{n+1}$ as it was done in
\cite{mish-noor}. Denote by $H_n$ the Hilbert space where the matrices
$\s_n(a)$, $\s_n(c)$ and $\s(b)$ act and put $H=\oplus_n H_n$ (for
negative $n$ put $\s_n(a)=\s_n(c)=I$). Let $F$ be a shift on $H$.
Then the matrices $\oplus_n\s_n(a)$, $\oplus_n\s_n(c)$
$\oplus_n\s(b)$ and $F$ generate a representation of $\G\times{\bf Z}$ into
${\cal Q}$ with necessary property.

\medskip
Remark that in our example the absence of the AGA property is related to
torsion. It would be interesting to know whether torsion-free groups
always have AGA.

\bigskip
{\bf Acknowledgement.}
The present paper was prepared with the partial
support of RBRF (grant N 96-01-00276) and INTAS.
I am grateful to A.~S.~Mishchenko for useful discussions.



\vspace{1cm}
{\small
\noindent
V.~M.~Manuilov\\
Dept. of Mech. and Math.,\\
Moscow State University,\\
Moscow, 119899, RUSSIA\\
e-mail: manuilov@mech.math.msu.su

}

\end{document}